\documentclass[10pt,a4paper,twoside,notitlepage]{article}

\usepackage{amsfonts}
\usepackage{amsmath}
\usepackage{amssymb}
\usepackage{color}
\usepackage{mathrsfs}

\newtheorem{thm}{Theorem}
\newtheorem{lem}{Lemma}[section]

\newtheorem{prop}{Proposition}[section]
\newtheorem{defn}{Definition}[section]
\newtheorem{rem}{Remark}[section]

\newcommand{\R}{\mathbb{R}}
\newcommand{\N}{\mathbb{N}}
\newcommand{\ve}{\varepsilon}
\newcommand{\vp}{\varphi}
\newcommand{\n}{\noindent}

\newcommand{\Div}{\mathrm{div}}
\newcommand{\Leb}{\mathcal{L}}
\newcommand{\con}{\mathbf{C}}
\newcommand{\Lsp}{L}
\newcommand{\prob}{\mathscr{P}}
\newcommand{\F}{{\frak F}}

\begin{document}

\title{First order mean field games in crowd dynamics}

\author{
  \vspace{1cm}
    {\scshape Fabio S. Priuli}
    \thanks{Istituto per le Applicazioni del Calcolo ``M. Picone'', C.N.R., Rome, Italy ({\tt f.priuli@iac.cnr.it})}
}

\date{\today}

\maketitle

\begin{abstract}
In this paper we study a two dimensional crowd model where pedestrian velocity consists of two elements: a non--local interaction term, modeling the effect of other walkers on each individual, and a control term. This latter term can be chosen by pedestrians so that their resulting path is optimal w.r.t. a suitable cost criterion. Under the assumption that pedestrians can forecast the effect of their choices on the evolution of the whole crowd, it is natural to consider a mean field system coupling the continuity equation, which describes the evolution of the density of pedestrians, with the HJB equation for the optimization problem. We show that such coupled system admits solution and we interpret the solution in terms of the pedestrian models. We also extend this results to the case where multiple populations of pedestrians are present in the environment, each with its own objectives.\\

\n{\bf Key words.} Crowd dynamics, mean field games, conservation laws, Hamilton-Jacobi\\

\n{\bf AMS subject classifications.} 35L65, 49N90, 49Q10, 91D10
\end{abstract}

\section{Introduction}

In recent years, the modeling of pedestrian motion has been subject to a wide attention for its practical relevance. The complex behavior of human crowds, however, poses a significant challenge and the available results only give a partial picture of the great variety of behaviors that can be observed in real situations.

From the mathematical modeling point of view, the research on pedestrian dynamics is relatively young, but it has already produced a considerable amount of literature. Many models, of very different nature, have been proposed (microscopic, mesoscopic, macroscopic, multiscale, cellular automata, discrete choice), and some of them in multiple variants as well (first-order, second-order, differential, non differential, local, non--local, with or without contact-avoidance features). We point the reader to the following reviews~\cite{bellomo2011,duives2013,helbing2001}, and to the books~\cite{kachroo2008,CPTbook}, for a general account and comparison of the existing models.

The main feature of pedestrian crowd models, compared to other physical models or to vehicular traffic flows, is that pedestrians are ``active particles''. Namely they are not passively subject to external force fields like inert matter. Rather they observe the surroundings and elaborate a strategy before moving. In many cases, they also have a specific target and want to reach it either in minimal time or while optimizing some other personal performance criterion. 

In the paper~\cite{CrPrTo}, we have introduced several different behaviors for pedestrians, differentiated on the basis of \emph{pedestrian rationality}, i.e. on the different amount of information available to pedestrians when they have to choose their strategy. 
Namely, we assumed that crowd dynamics obey the mass conservation principle, which implies that
the density of pedestrians $\rho$ has to satisfy the continuity equation
\begin{equation}\label{eq:cont_eq}
\partial_t\rho+\mathrm{div}(\rho v) = 0\,,
\qquad\qquad
\rho(0,\cdot)=\rho_0\,,
\end{equation}
where $t>0$, $x\in\R^2$, $v=v(t,x)$ is the macroscopic velocity of the crowd
and $\rho_0$ its initial distribution. The velocity field $v$ is explicitly chosen as the superposition of two different effects
\begin{equation}\label{eq:v}
v=v_b^*+v_i\,.
\end{equation}
The \emph{interaction velocity} $v_i$ is the component that pedestrians cannot control because it has to do with the effect produced by other walkers. In particular, it does not stem from choices made by individuals but rather from external conditions the latter undergo, and it depends from the distribution $\rho$ itself (we will denote this dependence by writing $v_i=v_i[\rho]$). The \emph{behavioral velocity} $v_b^*$ is, instead, the component that pedestrians can control more directly, i.e., the one which they can effectively act on to implement their behavioral strategy. The choice of the notation using an asterisk has been done so to convey the idea that the pedestrians try to choose it in what they perceive as an \emph{optimal way}, with respect to some personal criterion. 

Different degrees of ``rationality'' can be adopted by pedestrians in choosing $v^*_b$. For instance, they can choose their path only based on their knowledge of the environment, while completely ignoring what others do (\emph{basic behavior}, with no rationality); Or they can consider at each time $\tau$ only the distribution $\rho(\tau,\cdot)$ of others, reacting to other walkers' choices as if they were a new obstacle on their path (\emph{rational behavior}, with no predictive capabilities).

\n Here, we focus on the case of \emph{highly rational pedestrians}, in the terminology of~\cite{CrPrTo}. In this case, we assume that each pedestrian has a complete view of the distribution $\rho(\cdot,\cdot)$ at every time, and that he/she can forecast with precision the reaction of others to his/her decisions. This means that the evolution of the crowd has to be computed as a whole, because each individual affects the choices of the others. It also means that pedestrians will choose their control $v_b^*$ so to obtain a global (in time and space) optimum at the end of their evolution. We thus enter the field of \textit{differential games} and we can expect that pedestrians look for \textit{Nash equilibrium strategies} when choosing their control $v_b^*$. We recall that a Nash equilibrium is a situation where no pedestrian can find a better path if the others keep their choice, and the situation is somehow ``stuck'', since no one finds it convenient to change unilaterally his/her strategy. 

The assumption that each pedestrian have a complete understanding of the consequence of his/her choices for all later times, may sound restrictive when dealing with real situations, and therefore debatable. On the one hand, it is not uncommon to see pedestrians who can forecast to some extent the evolution of a crowd and choose a long--term optimal behavior: This is the case when pedestrians know very well the environment they are moving in, possibly because they pass through it every day, and when they can use such familiarity to predict under which conditions a change of strategy may be rewarding, e.g., making more convenient to move towards a farther target, because closer ones are going to be too crowded soon. On the other hand, in many situations crowds do not behave in this manner: If the environment is not well known or if there is some hurry in reaching the fixed target, so that decision has to be taken very rapidly, it sounds unrealistic to expect long term predictions from individuals.

However, following~\cite{CrPrTo}, here we would like to adopt a different point of view. Rational models like the one depicted above are not used to \textit{describe} but rather to \textit{control} reality. More specifically, the crowd evolution $(t,x)\mapsto \rho(t,x)$ described by highly rational pedestrians aims to be used as a target in a crowd control process, where the real behavior is steered towards the more rational one. Control of the crowd evolution can be performed either by adding some adjustable parameters in the continuity equation (e.g. a source term related to further density inflow in the environment), or by means of additional constraints to the admissible velocity fields $v_b^*$. Alternatively, one can decide to act on the environment to control the crowd, by adding (controlled) obstacles on pedestrians' path so to positively affect their natural behavior in the new environment (we refer to~\cite{CrPrTo} for a more detailed discussion about environmental optimization and for numerical tests).

In order to investigate the above control problems, it becomes important to have a robust mathematical framework ensuring that our target dynamics exist and have the properties we desire. This is the problem we start to tackle down in this paper, by proving that such highly rational dynamics exist, provided both the velocity field $v_i$ and the performance criterion that pedestrians want to optimize are smooth enough. Namely, we consider the case of dynamics taking place on a fixed time horizon $[0,T]$ and on the whole space $\R^2$, and we assume that pedestrians have access to the whole $\rho(\cdot,\cdot)$ when they have to choose their behavioral velocity field $(t,x)\in[0,T]\times\R^2\mapsto v_b^*(t,x)$. To fix the ideas, assume that pedestrians starting at time $t$ from the position $x$ want to minimize a common cost ${\cal C}(t, x,v_b(\cdot))$ along their path, i.e. along solutions\footnote{In this paper, solutions to $\dot y=f(t,y(t),u(t))$ in an interval $I=[\alpha,\beta]\subset\R$ are understood in Carath\'{e}odory sense, i.e., absolutely continuous functions on $I$ satisfying $y(t)= y(\alpha)+\int_\alpha^t f(\sigma,y(\sigma),u(\sigma))\,d\sigma$ for a.e.~$t\in\, I$.} on $[t,T]$ to
\begin{equation}\label{eq:dyn_real}
\dot y(s) = v_b(s)+v_i[\rho(s,\cdot)](y(s))\,,
\qquad\qquad
y(t)=x\,.
\end{equation}

\n By defining the value function for the problem as
\begin{equation}\label{eq:value}
\phi(t,x):=\inf_{v_b(\cdot)\in\mathcal A} \mathcal C(t,x,v_b(\cdot))\,,\qquad\qquad (t,x)\in[0,T]\times\R^2\,,
\end{equation}
where ${\cal A}$ is a suitable set of admissible measurable open--loop controls, it is well known in control theory that such a function can be retrieved as the viscosity solution to a suitable Hamilton--Jacobi--Bellman (HJB) equation. Namely, $\phi$ is the viscosity solution to
\begin{equation}\label{eq:gen_hb}
-\partial_t\phi(t,x)+H[\rho](t,x,D_x\phi) = 0\,,
\end{equation}
where $H[\rho]$ is a suitable Hamiltonian whose precise expression depends on the choice of ${\cal C}$. In turn, the knowledge of the value function  allows pedestrians to choose an optimal feedback control $v_b^*$. We thus end up with a system of PDEs~\eqref{eq:cont_eq}--\eqref{eq:gen_hb} in the unknowns $(\rho,\phi)$, which is fully coupled because the Hamiltonian depends on the density $\rho$ and the velocity fields $v$ in~\eqref{eq:cont_eq} depends on $\phi$ through $v_b^*$.

Mathematically, such coupled system~\eqref{eq:cont_eq}--\eqref{eq:gen_hb} is a first order mean field game, along the lines of~\cite{LL1,LL2,CardaNotes}. This should not be completely surprising: Solutions to mean field games are strictly related to the limit of Nash equilibrium strategies for games with $N$ agents, when $N\to\infty$, and thus they have found application in a large number of situations where a transport or a diffusion process is coupled with an optimization procedure, like in mathematical finance and in segregation models (see~\cite{GLL} for further applications). In the context of crowd dynamics, mean field games have been used by~\cite{Dogbe,LW,BuDFMaWo}, but always in cases where the individual behavior~\eqref{eq:dyn_real} evolves according to a stochastic process. This assumption has the effect to simplify some aspects of the analysis, since the continuity equation~\eqref{eq:cont_eq} is replaced by a parabolic PDE which admits in general smoother solution, but no analogous approach can be pursued in the deterministic situation we are studying here.

In this paper, we first extend the general existence result by Lasry \& Lions~\cite{LL1,LL2} for mean field games in $\R^d$ to the case of systems~\eqref{eq:cont_eq}--\eqref{eq:gen_hb} with a non--local velocity field appearing inside~\eqref{eq:v}, and thus in the Hamiltonian. Then, we apply such result to $2$d crowd models, obtaining the existence of highly rational dynamics and rigorously justifying the usage of such crowd behavior as desirable target for pedestrian optimization processes, as in~\cite{CrPrTo}. Finally, we also extend our results to the case of dynamics with multiple populations, each one with different (and possibly conflicting) goals, so to model situations where different groups of pedestrians interact with each other in the same environment (crowds crossing each others, travelers with different destinations, etc.). In such a case, the velocity field $v$ in~\eqref{eq:v} depends on the population because both $v_b^*$ and $v_i$ can depend on it: The difference in the behavioral velocity accounts for the different goals populations might have; The difference in the interaction velocity accounts for possible different reactions of individuals when interacting with their similar and/or with people of a different population (heterophobic or heterophilic behaviors). The result we prove is essentially the same as the previous one: if the vector fields and the cost functionals are smooth enough, then a solution exists also for multiple populations.

The paper is organized as follows. In Section~\ref{sec:theory_n_pedestrians} we present the result about mean field games and its application to pedestrian crowd models. In Section~\ref{sec:proof} we give the proof of Theorem~\ref{thm:mfg}. Finally, Section~\ref{sec:fine} discusses the extension to several populations and some open problems.

\section{Main results}\label{sec:theory_n_pedestrians}

The first result of this section deals with the existence of solutions to the mean field system
\begin{equation}\label{eq:mfg}
\left\{
\begin{array}{l}
-\partial_t u(t,x)+\displaystyle{1\over 2}\,|D_xu(t,x)|^2-D_xu(t,x)\cdot F(x,\mu(t))-G(x,\mu(t))=0\\
\partial_t \mu(t)+\Div\left(\mu(t)\,\big(-D_xu(t,x)+F(x,\mu(t))\big)\right)=0\phantom{\displaystyle{1\over 2}}\\
\mu(0)=\mu_0\,,~u(T,x)=\Psi(x,\mu(T))\phantom{\displaystyle{1\over 2}}
\end{array}
\right.
\end{equation}
for $(t,x)\in(0,T)\times \R^d$. Notice that the first equation is a backward Hamilton--Jacobi--Bellman (HJB) equation for an unknown function $u$, with Hamiltonian $H(t,x,p)={1\over 2}|p|^2-p\cdot F-G$, while the second equation is a measure--valued continuity equation for an unknown curve of measures $t\mapsto\mu(t)$. We search for solution pairs $(u,\mu)$ of~\eqref{eq:mfg} in the sense of the following definition. 

\begin{defn}\label{defn:sol} A pair $(u,\mu)\in W^{1,\infty}((0,T)\times\R^d)\times\con([0,T]\,;\prob_1(\R^d))$ is a \emph{solution pair} for system~\eqref{eq:mfg} if
$u$ is a viscosity solution for the first equation on $(0,T)\times\R^d$ with final condition $u(T,x)=\Psi(x,\mu(T))$ in $\R^d$, and if the curve $t\mapsto \mu(t)$ satisfies the initial condition $\mu(0)=\mu_0$ and it is a solution of the second equation in the sense of distributions, i.e. for all test functions
$\vp\in\con_c^\infty((0,T)\times\R^d)$ there holds
$$
\int_0^T\!\!\int_{\R^d} \Big[\partial_s \vp(s,x)+\langle D_x \vp(s,x),-D_xu(s, x)+F(x,\mu(s))\rangle\Big]\,d\mu(s)(x)\,ds=0\,.
$$
In the special case where all measures $\{\mu(s)\}_{s\in[0,T]}$ are absolutely continuous w.r.t. $\Leb^d$ and their densities, still denoted with $\mu$, give a map $(s,x)\mapsto \mu(s,x)\in\Lsp^1((0,T)\times\R^d)$, 
we will say that we have a solution pair $(u,\mu)\in W^{1,\infty}((0,T)\times\R^d)\times\Lsp^1((0,T)\times\R^d)$. 
\end{defn}

\n Concerning the coefficients and data of~\eqref{eq:mfg}, the standing assumptions are the following ones.

\begin{description}
\item{\bf (H)} The functions
$F\colon \R^d\times\prob_1(\R^d)\to\R^d$ and
$G,\Psi\colon \R^d\times\prob_1(\R^d)\to\R$ 
are continuous. Moreover, for any fixed probability measure $\mu\in\prob_1(\R^d)$, there hold
$F(\cdot,\mu)\in\con^2(\R^d\,;\R^d)$ and $G(\cdot,\mu),\Psi(\cdot,\mu)\in\con^2(\R^d\,;\R)$ with 
$$
\|F(\cdot,\mu)\|_{\con^2}\leq C\,,
\qquad \qquad
\|G(\cdot,\mu)\|_{\con^2}\leq C\,,
\qquad \qquad
\|\Psi(\cdot,\mu)\|_{\con^2}\leq C\,,
$$
for a constant $C$ independent from $\mu$; and for any fixed point $x\in\R^d$ the functions $F(x,\cdot)\colon \prob_1(\R^d)\to\R^d$ and $G(x,\cdot)\colon \prob_1(\R^d)\to\R$ are Lipschitz continuous with Lipschitz constant $C$ independent on the point $x$. Finally, the initial data $\mu_0$ is a measure such that $\mu_0\ll\Leb^d$, and its density still denoted by $\mu_0$ satisfies $\|\mu_0\|_\infty\leq C$ and $\mathrm{supp}~\mu_0\subseteq B(0,C)$.
\end{description}

\n Here and in the following, $\prob_1(\R^d)$ denotes the space of probability measures $\mu$ such that
$\int_{\R^d} |x-\bar x|\,d\mu(x)<+\infty$ for some $\bar x\in\R^d$, and it is considered as a metric space endowed with the $1$--Wasserstein metric, which can be defined as
follows
\begin{equation}\label{eq:wasserstein}
W_1(\rho,\sigma):=\sup\left\{\int_{\R^d}\!\!\vp\,d(\rho-\sigma)~;~\vp\colon\R^d\to\R\mbox{ is Lipschitz with }\mathrm{Lip}(\vp)\leq 1\right\}\,.
\end{equation}
Observe that assumptions {\bf (H)} are a bit more restrictive than it is really necessary, and they have been chosen for the sake of simplicity: One could replace the $\con^2$ regularity w.r.t. $x$ in {\bf (H)} with the weaker assumptions of Theorem~7.4.12 in~\cite{CS} and the following result would still be true.

\begin{thm}\label{thm:mfg} Assume that the functions $F$, $G$, $\Psi$ and the measure $\mu_0$ satisfy {\bf (H)}. Then, system~\eqref{eq:mfg} admits a solution pair 
$(u,\mu)\in W^{1,\infty}((0,T)\times\R^d)\times\Lsp^1((0,T)\times\R^d)$ in the sense of Definition~\ref{defn:sol}. In fact, $u$ is also semiconcave on $(0,T)\times\R^d$ and, for all $s\in[0,T]$, $\mu(s,\cdot)$ has compact support and satisfies $\mu(s,\cdot)\in\Lsp^\infty(\R^d)$.
\end{thm}

The proof of Theorem~\ref{thm:mfg} is deferred to Section~\ref{sec:proof}. Here we only anticipate that such proof follows the same pattern as the proofs by Lasry \& Lions~\cite{LL1,LL2}, using a fixed point argument on a compact map ${\cal T}$, defined on a suitable convex subset of $\con([0,T],\prob_1(\R^d))$. The most elaborate part of the construction is the proof that ${\cal T}$ is well--defined, because we need to prove existence and regularity of the solutions to the two equations in~\eqref{eq:mfg}.

In the rest of this section, instead, we apply Theorem~\ref{thm:mfg} to prove the existence of solutions for a particular two--dimensional, macroscopic model for crowd dynamics. We consider as domain where pedestrians move (the so--called \textit{walking area}) the whole $\R^2$, and we assume that the evolution takes place on a fixed time horizon $[0,T]$. We denote the time variable by $t$, the space variable by $x$, and the spatial density of the continuum representing the crowd by $\varrho=\varrho(t,x)$. Pedestrian dynamics are assumed to obey the mass conservation principle, which implies that $\varrho$ is a solution to the continuity equation~\eqref{eq:cont_eq}, with initial distribution $\varrho_0$, which we assume to be in $\Lsp^1(\R^2)\cap\Lsp^\infty(\R^2)$ and to have compact support, and with velocity field $v=v(t,x)\in\R^2$ given by~\eqref{eq:v}, which prescribes the macroscopic speed of the crowd. In particular, owing to the mass conservation principle, we have that $\|\varrho(t,\cdot)\|_{\Lsp^1(\R^2)}\equiv\|\varrho_0\|_{\Lsp^1(\R^2)}$ for all $t\in[0,T]$. Thus, we associate to $\varrho(t,\cdot)$ an absolutely continuous probability measure 
$$
\rho(t):={\varrho(t,\cdot)\over \|\varrho_0\|_{\Lsp^1(\R^2)}}\,\Leb^2\,,
$$ 
and we obtain that the corresponding curve of measures $t\mapsto\rho(t)$ is a weak solution of the measure--valued continuity equation~\eqref{eq:cont_eq}, still driven by the velocity field $v$ and with $\rho_0=\varrho_0\,\Leb^2/ \|\varrho_0\|_{\Lsp^1}$.

Since we want to model pedestrians which elaborate actively a behavioral strategy for their movements, we refrain from invoking inertial principles for obtaining their velocity $v$ (which would lead to a second order model). We prefer to adopt a first order modeling approach, in which pedestrian velocity is genuinely modeled by taking into account the superposition of two basic contributions: On the one hand, the desire to reach a specific destination (e.g., an exit point) in a suitably \emph{optimal way}, which generates a \textit{behavioral velocity} $v_b^*$; On the other hand, the necessity to avoid collisions with other surrounding walkers, which generates an \textit{interaction velocity} $v_i$. This results in the expression~\eqref{eq:v}, $v=v_b^*+v_i$.

Observe that the two summands in $v$ have a completely different interpretation. We start from the interaction velocity $v_i$. Here, following~\cite{CPT10,CPT11,TosinFrasca,PiccoliRossi,CrPrTo}, we consider a \emph{non--local} interaction velocity which depends on the distribution of pedestrians in the walking area, so to describe short- and medium-range influence of nearby walkers on a pedestrian. 
Namely, if a distribution of pedestrians on $\R^2$ is given, as a measure $\sigma\in\prob_1(\R^2)$, we assume that the interaction velocity is a vector field $v_i[\sigma]\colon\R^2\to\R^2$ whose pointwise value in $x$ depends on the $\sigma$--measure of a suitable neighborhood of $x$. To fix the idea, you can think to an integral vector field as $v_i[\sigma](x)=\int_{\R^2}\tau(x,y)\,d\sigma(y)$, for a suitable family of $\sigma$--summable vector fields $\tau(x,\cdot)$ in $\R^2$. For our purposes, the precise assumptions on the velocity field $v_i$ are the following ones.
\begin{description}
\item{\bf (H1)} The function $v_i\colon \prob_1(\R^2)\to \con^2(\R^2\,;\R^2)\cap\Lsp^\infty(\R^2\,;\R^2)$ is continuous as a map between the metric spaces $\prob_1(\R^2)$, endowed with the $1$--Wasserstein metric~\eqref{eq:wasserstein}, and $\con(\R^2\,;\R^2)$. 
Moreover, there exists a suitable constant $C'> 0$ such that the following holds: 
For every measure $\sigma$, the resulting map $v_i[\sigma](\cdot)$ is in $\con^2(\R^2\,;\R^2)$ and $\|v_i[\sigma]\|_{\con^2}\leq C'$;
For every $x\in\R^d$, the resulting map $v_i[\cdot](x)$ satisfies $|v_i[\sigma](x)-v_i[\sigma'](x)|\leq C'\,W_1(\sigma,\sigma')$ for every $\sigma,\sigma'\in\prob_1(\R^2)$.
\end{description}

Concerning the behavioral velocity $v_b^*$, we treat it as a control function to be chosen by pedestrians to reach their objective in what they perceive as an optimal way. In other words, $v_b^*$ is a control strategy and we have to describe in which way individuals select it during the evolution. As mentioned in the Introduction, here we assume that pedestrians are \emph{highly rational}, i.e. that they are capable to predict exactly the evolution of the density measure $\rho(t)$ during the whole time horizon $[0,T]$ and that they can forecast the long term effect of their choices on the other pedestrians too. 
Moreover, we assume that pedestrians want to minimize the following performance criterion
\begin{equation}\label{eq:cost2}
{\cal C}(t, x,v_b(\cdot))~:=~\int_t^T \left({1\over 2} |v_b(s)|^2+\ell[\rho(s)](y(s))\right)\,ds+h[\rho(T)](y(T))\,,
\end{equation}
where, for a fixed curve $t\mapsto\rho(t)$ with values in $\prob_1(\R^2)$, $\ell$ is a positive running cost, accounting for the cost of a specific path, $h$ is a bounded below exit cost, accounting for the possible failure to reach the target at time $T$, and $y(\cdot)$ is the microscopic path followed by the pedestrian, i.e. the solution to
\begin{equation}\label{eq:dyn_perc}
\dot y(t) = v_b(t)+v_i[\rho(t)](x)\,.
\end{equation}
The dependence of the cost functions $\ell,h$ from the distribution of pedestrians $\rho=\varrho\Leb^2$ is again non--local and it can represent a way to penalize aggregation. In the following, we assume these hypotheses are verified:
\begin{description}
\item{\bf (H2)} The functions $\ell,h\colon \prob_1(\R^2)\to \con^2(\R^2\,;\R)\cap\Lsp^\infty(\R^2\,;\R)$ are continuous as maps between the metric spaces $\prob_1(\R^2)$, endowed with the $1$--Wasserstein metric~\eqref{eq:wasserstein}, and $\con(\R^2\,;\R)$.
Moreover, there exists suitable constants $C'> 0$ and $\ve>0$ such that the following holds: 
For every measure $\sigma$, the resulting maps $\ell[\sigma]$ and $h[\sigma]$ are in $\con^2(\R^2\,;\R)$ with 
$$
\ell[\sigma]\geq\ve>0\,,
\qquad\qquad
\|\ell[\sigma]\|_{\con^2}\leq C'\,,
\qquad\qquad
\|h[\sigma]\|_{\con^2}\leq C'\,;
$$
For every $x\in\R^d$, the resulting map $\ell[\cdot](x)$ satisfies for every $\sigma,\sigma'\in\prob_1(\R^2)$
$$
|\ell[\sigma](x)-\ell[\sigma'](x)|\leq C'\,W_1(\sigma,\sigma')\,.
$$ 
\end{description}

Once pedestrians are aware of either $\varrho(\cdot,\cdot)$ on $[0,T]\times\R^2$ or $\rho(\cdot)$ on $[0,T]$, the actual construction of the control strategy $(t,x)\in[0,T]\times\R^2\mapsto v_b^*(t,x)$, in feedback form, can proceed as follows. Each individual defines his/her value function $\phi$ for the problem~\eqref{eq:dyn_perc}--\eqref{eq:cost2} as in~\eqref{eq:value}, and find $\phi$ by solving an HJB equation like~\eqref{eq:gen_hb} with Hamiltonian
\begin{equation}\label{eq:hamilton}
H(t,x,D_x\phi):=\max_{v_b \in \R^2}\left\{ -(v_b+v_i[\rho(t)](x)\big)\cdot D_x\phi(t,x)-{1\over 2} |v_b|^2-\ell[\rho(t)](x)\right\}\,,
\end{equation}
and final condition $\phi(T,x)=h[\rho(T)](x)$. Once $\phi$ has been found, owing to Pontryagin minimum principle, the optimal feedback control $v_b^*(t,x)$ can be selected by choosing an element which realizes the maximum in~\eqref{eq:hamilton}. Given the particular form of the cost considered here, the expression can be actually given explicitly as
$$
v_b^*(t,x)=-D_x\phi(t,x)\,,
$$
and thus we are facing the following \emph{forward--backward} system of PDEs 
\begin{equation}\label{eq:mfg2}
\left\{
\begin{array}{l}
-\partial_t\phi(t,x)+\displaystyle\,{1\over 2}\,|D_x\phi(t,x)|^2-v_i[\rho(t)](x)\cdot D_x\phi(t,x)-\ell[\rho(t)](x)= 0\,,\\
\partial_t \rho(t)+\Div\Big(\rho(t)\,(-D_x\phi(t,x)+v_i[\rho(t)](x))\Big)=0\,,\phantom{\displaystyle{1\over 2}}
\end{array}
\right.
\end{equation}
with initial condition $\rho(0,\cdot)=\rho_0$ and final condition $\phi(T,\cdot)=h[\rho(T)]$ on $\R^2$. Note that the system is of the form~\eqref{eq:mfg} and that the equations are \emph{fully coupled} because of the presence of the whole density $\rho$ in the first equation as well as the presence of $\phi$ in the second one. Our result for this crowd model is the following one.

\begin{thm}\label{thm:crowd} Assume that $v_i$ satisfies {\bf (H1)}, that $\ell$ and $h$ satisfy {\bf (H2)} and that $\varrho_0\in\Lsp^1(\R^2)\cap\Lsp^\infty(\R^2)$ has compact support.
Then, system~\eqref{eq:mfg2} admits a solution pair $(\phi,\varrho)\in
W^{1,\infty}((0,T)\times\R^2)\times\Lsp^1((0,T)\times\R^2)$ in the sense of Definition~\ref{defn:sol}. In fact, $\phi$ is also semiconcave on $(0,T)\times\R^2$ and, for all $s\in[0,T]$, $\varrho(s,\cdot)$ has compact support and belongs to $\Lsp^1(\R^2)\cap\Lsp^\infty(\R^2)$.
\end{thm}

\smallskip

\n{\it Proof.} Assumptions {\bf (H1)} and {\bf (H2)} ensure that the functions $F(x,\sigma):=v_i[\sigma](x)$, $G(x,\sigma):=\ell[\sigma](x)$ and $\Psi(x,\sigma):=h[\sigma](x)$, defined for all $(x,\sigma)\in\R^2\times\prob_1(\R^2)$, satisfy assumptions {\bf (H)}. Thus, we can apply Theorem~\ref{thm:mfg} to system~\eqref{eq:mfg2} to obtain the conclusion.~~$\diamond$

\smallskip

In view of Theorem~\ref{thm:crowd}, we can offer a justification for the choice of ``highly rational'' to describe this pedestrian behavior (cf also~\cite{CrPrTo}). Indeed, when a solution pair $(\phi, \varrho)$ of~\eqref{eq:mfg2} exists, such a pair it is likely to represent the crowd distribution and the corresponding performance outcome for the cost ${\cal C}$ when the pedestrians are close to a Nash equilibrium between themselves. Indeed, it is believed (and rigorously proved in several situations) that solutions of the mean field limit system~\eqref{eq:mfg2} give \emph{approximate Nash equilibria} for the crowd when the number of people in the walking area is large enough. More precisely, if we consider the control problem~\eqref{eq:dyn_perc}--\eqref{eq:cost2} as a game among $N$ competing microscopic agents (viz. pedestrians), indistinguishable, distributed according to the measure $\varrho$, and trying to minimize their common cost ${\cal C}$, then we expect that $\phi$ represents the outcome of each pedestrian at time $T$ and that no one of them can improve it more than $O(1/N)$ by means of unilateral changes of strategy. 
 
\begin{rem}
Observe that the assumptions {\bf (H1)} on $v_i$ are quite restrictive, because they require uniform Lipschitz continuity w.r.t. the measure and because they rule out interaction velocities which depend on the pointwise values of the density $\varrho$. However, such assumptions are not too restrictive for the application to pedestrian flows we have in mind: The explicit expressions of $v_i$ used in~\cite{CPT11,CrPrTo} do satisfy {\bf (H1)}, for suitable choice of the repulsion functions. The usual expression for $v_i[\varrho\Leb^2]$, in the cited papers, is given by
\begin{equation}\label{eq:vi}
v_i[\varrho\Leb^2](x)=\int_{B(x,R)} \F(y-x)\varrho(y)\,dy\,,
\end{equation}
for some fixed $R>0$, so to model the ability of pedestrians to scan their surroundings and react to the distribution of other people there, and for a function $\F:\R^2\to\R^2$ which models the reaction of an individual in $x$ to another individual in $y$. For collision avoidance purposes, $\F$ is typically considered as a repulsion and taken as $\F(\xi)=-{\kappa \over |\xi|^2}\,\xi$, for a suitable constant $\kappa>0$ tuning the strength of the repulsion. Notice that, since the distance between two pedestrians cannot decrease below a certain threshold $r_o\in(0,R)$ (due to the physical size of pedestrian bodies), it is always possible to cut--off $\F$ when $|\xi|<r_o$ and to mollify $\F$ with a sufficiently smooth function $\eta_R\colon\R^2\to [0,\,1]$, compactly supported in $B(0,R)$ but identically equal to $1$ in $B(0,R-\epsilon)$, for some $\epsilon\in(0,r_0)$. In this way, the singularity at $\xi = 0$ can be avoided and any desired regularity can be imposed to $\F$, without affecting the validity of the model (cf~\cite{cristiani2012CDC}). In particular, its global Lipschitz continuity implies the corresponding property for the map $\mu\mapsto v_i[\mu](x)$ independently by $x\in\R^d$.\\
We conclude by observing that in~\eqref{eq:vi} no anisotropy of the interactions is accounted, so that pedestrians are affected in their choices by both people at their front and people at their back. This is a restriction, compared to~\cite{CPT11,CrPrTo}, that we hope to relax in a future work.
\end{rem}

\begin{rem}
The non--local interaction velocity we consider here is fundamentally different from the ones considered in a large part of the literature for pedestrian crowds (e.g.~\cite{BuDFMaWo,CHM,HB,HB2,Hughes,LW}). Indeed, many authors assume that the macroscopic velocity $v(t,x)$ is given by an expression of the form $s(\varrho) w(t,x)$, where $w(t,x)$ is a prescribed unit vector field accounting for the target of pedestrians, and $s(\varrho)$ is a scalar factor that decreases when the pointwise density $\varrho(t,x)$ increases, so to slow down pedestrians when the environment is crowded. The actual expression of the function $s(\varrho)$ is what characterizes each model (a popular choice is to fix suitable positive constants $V_M$ and $\varrho_M$ and take $s(\varrho)=V_M(1-\varrho/\varrho_M)$), and it plays an entirely analogous role to the one of fundamental diagrams in $1$d models for vehicular traffic flows. However, in the context of $2$d pedestrian flow models any expression of $s(\varrho)$ has yet to find a rigorous justification, from the kinematic point of view. 

\n Therefore, we preferred to follow the approach initiated by~\cite{CPT10}, using an \emph{intrinsically multidimensional model} for interactions, where individuals are not only affected by the people in the same position or along the same path, but by all people present in a $2$d area surrounding him/her.
\end{rem}

\section{Proof of Theorem~\ref{thm:mfg}}\label{sec:proof}

As mentioned in Section~\ref{sec:theory_n_pedestrians}, the proof proceeds along the lines of~\cite{LL1,LL2} and makes use of many standard results from control theory. Therefore, we skip some of the details (giving precise references to~\cite{CS,CardaNotes,LL1,LL2} where necessary), and we focus our attention on the parts which require some specific treatment, due to the presence of the nonlocal term $F$ in the Hamiltonian.

\n In subsection~\ref{sec:HJB}, we study the HJB equation for a fixed curve of measures, proving that it has a unique semiconcave solution $u$ and analyzing a related optimization problem for which $u$ is the value function. Then, in subsection~\ref{sec:CLAW} we deal with the continuity equation, proving that it admits a unique weak solution, consisting of the push--forward of $\mu_0$ along the integral curves of $-D_xu+F$. Finally, subsection~\ref{sec:PROOF} contains the fixed point argument needed to conclude.

\subsection{The Hamilton--Jacobi--Bellman equation in~\eqref{eq:mfg}}\label{sec:HJB}

We fix a Lipschitz continuous curve of probability measures $t\mapsto m(t)\in\con([0,T]\,;\prob_1(\R^d))$ such that $m(0)=\mu_0$, and we focus our attention to the corresponding first equation in~\eqref{eq:mfg} with its final condition. By setting for all $(t,x)\in[0,T]\times\R^d$
\begin{equation}\label{eq:fixed_funct}
f(t,x):= F(x,m(t))\,,
\qquad\qquad
g(t,x):= G(x,m(t))\,,
\qquad\qquad
\psi(x):= \Psi(x,m(T))\,,
\end{equation}
we can rewrite such equation as
\begin{equation}\label{eq:hjb}
\left\{
\begin{array}{l}
-\partial_t u(t,x)+\displaystyle\,{1\over 2}\,|D_xu(t,x)|^2-D_xu(t,x)\cdot f(t,x)-g(t,x)=0\,,\\
u(T,x)=\psi(x)\,.\phantom{\displaystyle{1\over 2}}
\end{array}
\right.
\end{equation}
We can interpret this equation as the HJB equation for the value function of a suitable optimization problem: Consider the nonlinear control system
\begin{equation}\label{eq:control}
\dot x(s) = \alpha(s) + f(s,x(s))\,,
\qquad\qquad
x(t)=x\,,
\end{equation}
with cost
\begin{equation}\label{eq:cost}
{\cal C}(t,x,\alpha(\cdot)):= \int_t^T \left({1\over 2} |\alpha(s)|^2+g(s,x(s))\right)\,ds+\psi(x(T))
\end{equation}
where $x(\cdot)$ is the unique Carath\'eodory solution to~\eqref{eq:control}, corresponding to the control $\alpha\in\Lsp^2([0,T]\,;\R^d)$.
By defining the value function as
\begin{equation}\label{eq:value_funct}
u(t,x):=\inf_{\alpha\in\Lsp^2([0,T]\,;\R^d)}\left\{\int_t^T \left({1\over 2} |\alpha(s)|^2+g(s,x(s))\right)\,ds+\psi(x(T))\right\}\,,
\end{equation}
we can use standard arguments from control theory (with an unbounded set of controls) to prove that $u$ is the unique bounded uniformly continuous viscosity solution to the (local) equation~\eqref{eq:hjb} on $(0,T)\times\R^d$. In fact, assumptions {\bf (H)} on $F,G,\Psi$ yield that $f(t,\cdot)$, $g(t,\cdot)$ and $\psi$ in~\eqref{eq:fixed_funct} are functions in $\con^2(\R^d\,;\R^d)$ with 
\begin{equation}\label{eq:H'}
\|f(t,\cdot)\|_{\con^2}\leq C\,,
\qquad\qquad
\|g(t,\cdot)\|_{\con^2}\leq C\,,
\qquad\qquad
\|\psi\|_{\con^2}\leq C\,,
\end{equation}
where $C$ is the constant, independent from $t$, appearing in {\bf (H)}. Thanks to~\eqref{eq:H'} and to the particular structure of dynamics and cost, we can prove that optimal controls exist for~\eqref{eq:control}--\eqref{eq:cost} (see e.g. Corollary~7.4.7 in~\cite{CS}). Namely, for every $(t,x)\in [0,T]\times\R^d$, there exists $\alpha^*\in\Lsp^2([0,T]\,;\R^d)$ such that
\begin{equation}\label{eq:inf_attained}
{\cal C}(t,x,\alpha^*(\cdot))=\inf_{\alpha\in\Lsp^2([0,T]\,;\R^d)} {\cal C}(t,x,\alpha(\cdot))\,,
\end{equation}
so that the infimum in~\eqref{eq:value_funct} is in fact a minimum. Moreover, any optimal open--loop control must satisfy the following necessary conditions, which follows from Pontryagin minimum principle (cf Theorem~7.4.17 in~\cite{CS}). 

\begin{prop}\label{prop:pmp} Let $\alpha^*\in\Lsp^2([t,T]\,;\R^d)$ be an optimal control for the problem~\eqref{eq:control}--\eqref{eq:cost}, let $x^*\colon[t,T]\to\R^d$ be the corresponding optimal trajectory and assume that $f,g$ satisfy~\eqref{eq:H'}. Then, $\alpha^*$ belongs to $\con^1([t,T]\,;\R^d)$ and it is a classical solution to the linear system
\begin{equation}\label{eq:opt_ctrl}
\dot \alpha^*(s) = D_xg(s,x^*(s)) -\,\alpha^*(s)\cdot D_xf(s,x^*(s))\,,
\qquad\qquad
\alpha^*(T)=-\nabla\psi(x^*(T))\,.
\end{equation}
In particular, there exists a constant $C_1>0$, only depending on the constant $C$ in~\eqref{eq:H'}, such that $\|\alpha^*\|_{\Lsp^\infty([t,T]\,;\R^d)}\leq C_1$.
\end{prop}

\smallskip

\n{\it Proof.} We recast the minimization problem~\eqref{eq:control}--\eqref{eq:cost} in $\R^{d+1}$ by introducing an additional variables $x_{d+1}$ which satisfies
$\dot x_{d+1}(s)=\,{1\over 2} |\alpha(s)|^2+g(s,x(s))$ and $x_{d+1}(t)=0$.
Then, we deal with an optimization problem for the variable $z:= (x,x_{n+1})\in\R^{d+1}$ with dynamics and cost criterion respectively given by
$$
\dot z(s)=
\!\left(
\begin{array}{c}
\alpha_1(s)+f^1(s,x(s))\\
\vdots\\
\alpha_d(s)+f^d(s,x(s))\\
{1\over 2} |\alpha(s)|^2+g(s,x(s))
\end{array}
\right)\!=: \frak{G}(s,z(s),\alpha(s))\,,
\quad
\frak{H}(z(T)):= x_{d+1}(T)+\psi(x(T)) \,.
$$
Owing to Pontryagin minimum principle, for any optimal pair $(\alpha^*,z^*)$, there exists an absolutely continuous arc $q=(p,p_{d+1})\colon[0,T]\to\R^{d+1}$ which solves the adjoint system
\begin{equation}\label{eq:adjoint}
\dot q(s) = -\,q(s)\cdot D_z\frak{G}(s,z^*(s),\alpha^*(s))\,,
\qquad\qquad
q(T)=D_z\frak{H}(z^*(T))
\end{equation}
and satisfies for almost every $s\in [t,T]$ the minimality condition
\begin{equation}\label{eq:optimal}
q(s)\cdot \frak{G}(s,z^*(s),\alpha^*(s))= \min \big\{q(s)\cdot \frak{G}(s,z^*(s),\omega)~;~\omega\in\R^d\big\}\,.
\end{equation}
Noticing that the extended dynamic $\frak{G}(s,z^*(s),\alpha^*(s))$ does not depend on $x_{d+1}$, it is immediate to deduce from~\eqref{eq:adjoint} that $p_{d+1}\equiv 1$.
Then, condition~\eqref{eq:optimal} reduces to
 ${1\over 2}|\alpha^*(s)|^2+p(s)\cdot\alpha^*(s)=\min \big\{{1\over 2}| \omega|^2+p(s)\cdot \omega
~;~\omega\in\R^d\big\}$,
which implies $\alpha^*(s)=-p(s)$ for $\Leb^1$--a.e. $s\in [t,T]$. 
Since $(\alpha^*,z^*)=(\alpha^*,x^*,x_{d+1}^*)$ is optimal for the extended problem if and only if $(\alpha^*,x^*)$ is optimal for the original one, with value ${\cal C}(t,x,\alpha^*)=x_{d+1}^*$, the conclusion follows by~\eqref{eq:adjoint}. 
The estimate on $\|\alpha^*\|_{\Lsp^\infty}$ follows from~\eqref{eq:H'} and Gronwall's Lemma.~~$\diamond$

\smallskip

In turn, boundedness of the optimal controls allows to prove that $u$ is Lipschitz continuous w.r.t. $(t,x)$ on $(0,T)\times\R^d$ and semiconcave w.r.t. $x$ on $\R^d$, for any time $t\in(0,T)$, with a constant $C_2\geq0$ which only depends on $C$ in~\eqref{eq:H'}, but neither on the functions $f,g,\psi$ nor on the time $t$ (cf Theorem~7.4.12 in~\cite{CS}). Semiconcavity of $u$ w.r.t. $x$ means that for all $x,h\in\R^d$ there holds
\begin{equation}\label{eq:semiconc}
u(t,x+h)+u(t,x-h)-2u(t,x)\leq C_2\, |h|^2\,.
\end{equation}
This uniform semiconcavity property for the value function~\eqref{eq:value_funct} will play an important role in Section~\ref{sec:PROOF}. 
Observe also that, by combining the (uniform) Lipschitz continuity of $F(x,\cdot), G(x,\cdot)$ with constant $C$ and the Lipschitz continuity of the curve $t\mapsto m(t)$ with constant $L_m$, we obtain for every $x\in\R^d$ 
$$
|f(t,x)-f(s,x)|\leq C\,L_m\,|t-s|\,,
\qquad\qquad
|g(t,x)-g(s,x)|\leq C\,L_m\,|t-s|\,,
$$ 
which yield that in fact $u$ is semiconcave w.r.t. both variables in $(0,T)\times\R^d$, this time with a constant which depends on both $C$ in~\eqref{eq:H'} and $L_m$ (the proof follows again by the standard arguments used for Theorem~7.4.12 in~\cite{CS}).
We collect in the following proposition some regularity results for semiconcave functions, whose proof can be found in~\cite{CS}.

\begin{prop}\label{prop:semiconc_lip} Let $v\colon (0,T)\times\R^d\to\R$ be a semiconcave function. Then, $v$ is locally Lipschitz continuous and, hence, it is differentiable almost everywhere in $[0,T]\times\R^d$.\\ 
Moreover, for every point $z=(t,x)\in(0,T)\times\R^d$ the \emph{superdifferential} of $v$ at $z$
\begin{equation}\label{eq:superdiff}
D_{t,x}^+v(z):= \left\{p\in\R^{d+1}~;~\limsup_{\zeta\to z}\,{v(\zeta)-v(z)-\langle p, \zeta-z\rangle\over |\zeta-z|}\,\leq 0\right\}
\end{equation}
is nonempty, compact and convex and it is the closed convex hull of the set $D_{t,x}^*v(z)$ of \emph{reachable gradients} of $v$ at $z$, i.e.
\begin{equation}\label{eq:reach_grad}
D_{t,x}^*v(z):= \left\{p\in\R^{d+1}~;~\exists\,z_n\to z \mbox{ such that $D_{t,x}v$ exists at $z_n$ and } D_{t,x}v(z_n)\to p\right\}\,.
\end{equation}
Thus, $v$ is differentiable at $(t,x)$ if and only if $D_{t,x}^+v(t,x)$ reduces to a singleton and in this case $D_{t,x}^+v(t,x)=D_{t,x}^*v(t,x)=\{D_{t,x}v(t,x)\}$.
Finally, if $(v_n)$ is a sequence of uniformly semiconcave functions with pointwise limit $v$, then $v$ is semiconcave, $v_n\to v$ locally uniformly, and $D_{t,x}v_n(t,x)\to D_{t,x}v(t,x)$ for a.e. $(0,T)\times\R^d$.
\end{prop}

\medskip

\n The same properties hold, of course, for the function $v(t,\cdot)$ seen as semiconcave function on $\R^d$.

Next lemma sums up the main properties of the set of optimal controls
\begin{equation}\label{eq:opt_set}
\mathscr{A}(t,x):=\left\{\alpha\in\Lsp^2([0,T]\,;\R^d)~;~\alpha\mbox{ is optimal for the problem~\eqref{eq:control}--\eqref{eq:cost}}\right\}\,,
\end{equation}
and its connection with the $D^+_xu$. Recall that $\mathscr{A}(t,x)\neq\emptyset$ because it contains $\alpha^*$ from~\eqref{eq:inf_attained}.

\begin{lem}\label{lem:opt_set} Let $(t,x)\in[0,T]\times\R^d$. Then, the following facts hold.
\begin{description}
\item{\it (i)} If $(t_n,x_n)\to(t,x)$ and $\alpha_n\in\mathscr{A}(t_n,x_n)$ satisfies $\alpha_n\rightharpoonup \alpha$ in $\Lsp^2([0,T]\,;\R^d)$, then $\alpha\in\mathscr{A}(t,x)$.
\item{\it (ii)} $D_xu(t,x)$ exists if and only if $\mathscr{A}(t,x)$ is a singleton $\{\alpha(\cdot)\}$ and if this is the case, then $D_xu(t,x)=-\alpha(t)$.
\item{\it (iii)} If $\alpha\in\mathscr{A}(t,x)$ and $x^*(\cdot)$ is the corresponding optimal solution to~\eqref{eq:control}, then for any $s\in(t,T)$ the set $\mathscr{A}(s,x^*(s))$ reduces to a singleton given by the restriction of $\alpha$ to $[s,T]$.
\end{description}
\end{lem}

\smallskip

\n{\it Proof.} {\it (i)} The optimization problem~\eqref{eq:control}--\eqref{eq:cost} with initial condition $x(t_n)=x_n$ $\Gamma$--converges to the problem with initial data $x(t)=x$ (see~\cite{BuFr}). This implies that the limit control $\alpha(\cdot)$ is optimal for the limit problem.

\n{\it (ii)} Assume that $D_xu(t,x)$ exists, let $\alpha(\cdot)\in\mathscr{A}(t,x)$ and $x(\cdot)$ be the corresponding trajectory with initial data $x(t)=x$. By using an argument analogous to the proof of Theorem~7.4.17 in~\cite{CS}, one verifies that the co--state $p(\cdot)$ corresponding to the optimal pair $(\alpha(\cdot), x(\cdot))$ satisfies $p(t)\in D^+_xu(t,x)$. But we have already seen that $p(\cdot)=-\alpha(\cdot)$, so it must be $-\alpha(t)\in D^+_xu(t,x)=\{D_xu(t,x)\}$. This implies that any optimal control would be a solution to the linear ODE~\eqref{eq:opt_ctrl} satisfying $\alpha(t)=-D_xu(t,x)$, so that $\mathscr{A}(t,x)$ must be singleton.

\n Viceversa, assume that $\mathscr{A}(t,x)$ is a singleton.
Fixed any reachable gradient $q\in D_x^*u(t,x)$, defined analogously to~\eqref{eq:reach_grad}, we have a sequence of points $x_n\to x$ such that $u(t,\cdot)$ is differentiable at $x_n$ and $D_xu(t,x_n)\to q$. From the previous implication, we know that $\mathscr{A}(t,x_n)=\{\alpha_n(\cdot)\}$ for all $n$ and, hence, the optimal pair $(\alpha_n(\cdot),x_n(\cdot))$ is the unique solution to the system
$$
\dot x_n(s) = \alpha_n(s) + f(s, x_n(s))\,,
\qquad\qquad
\dot \alpha_n(s) = D_xg(s,x_n(s))-\alpha_n(s) D_xf(s, x_n(s))\,,
$$
with initial conditions
$x_n(t) = x_n$ and $\alpha_n(t)= -D_x u(t,x_n)$. By letting $n\to\infty$, it is now clear that the sequence $(\alpha_n(\cdot),x_n(\cdot))$ converges uniformly to some pair $(\alpha(\cdot),x(\cdot))$ 
which solves the limit system
$$
\dot x(s) = \alpha(s) + f(s, x(s))\,,
\qquad\qquad
\dot \alpha(s) = D_xg(s,x(s))-\alpha(s) D_xf(s, x(s))\,,
$$
with conditions
$x(t) = x$ and $\alpha(t)= -q$.
Thus, $\alpha(\cdot)$ is optimal for $u(t,x)$ by part {\it (i)}. Being $\mathscr{A}(t,x)$ a singleton by assumption, we can conclude that $D_x^*u(t,x)$ must be a singleton too, and thus so is $D^+_xu(t,x)$, whence differentiability follows.

\n{\it (iii)} The dynamic programming principle implies that the restriction of $\alpha$ to $[s,T]$ is optimal for the problem with initial data $x(s)=x^*(s)$, so we just have to prove that $\mathscr{A}(s,x^*(s))$ reduces to a singleton. Thanks to {\it (ii)}, it is then enough to prove that $D_x u(s,x^*(s))$ exists for every $s\in(t,T)$. Given that $u(t,x)$ is semiconcave on $[0,T]\times\R^d$, we can equivalently prove that $D^+_xu(s,x^*(s))$ is a singleton, and this can be accomplished as in Theorem~7.3.16 in~\cite{CS}, by combining the dynamic programming principle with the strict convexity of the Hamiltonian $\,{1\over 2}\,|p|^2-p\cdot f-g$.~~$\diamond$

\smallskip

\n We conclude with a characterization of optimal controls, analogous to Lemma~4.11 in~\cite{CardaNotes}.

\begin{prop}\label{prop:iff_optimal_control} Let $(t,x)\in[0,T[\,\times\R^d$ and let $u(t,x)$ be given by~\eqref{eq:value_funct}. If $\alpha\in\mathscr{A}(t,x)$, then the associated optimal trajectory $x^*$ is a Carath\'eodory solution to the Cauchy problem
\begin{equation}\label{eq:opt_control}
\dot x(s) = -D_xu(s,x(s))+f(s,x(s))\,,
\qquad\qquad
x(t)=x\,.
\end{equation}
Conversely, if $x(\cdot)$ is a Carath\'eodory solution to~\eqref{eq:opt_control} defined on $[t,T]$, then the control function $\alpha(s):= \dot x(s)-f(s,x(s))$ is optimal for the problem~\eqref{eq:control}--\eqref{eq:cost}.
\end{prop}

\smallskip

\n{\it Proof.} From parts {\it (ii)} and {\it (iii)} of Lemma~\ref{lem:opt_set} we deduce that $D_xu(s,x^*(s))$ exists for all $s\in(t,T)$ and that $D_xu(s,x^*(s))=-\alpha(s)$. The conclusion of the first part is then immediate. \\
We omit the proof of the converse implication, because it is completely analogous to the one of Lemma~4.11 in~\cite{CardaNotes}. Indeed, the property only depends on the regularity of the Hamiltonian and not on its specific expression.~~$\diamond$

\subsection{The continuity equation in~\eqref{eq:mfg}}\label{sec:CLAW}

We now pass to the second equation of~\eqref{eq:mfg} with its initial condition. As in the previous section, let $t\mapsto m(t)\in\con([0,T]\,;\prob_1(\R^d))$ be a Lipschitz continuous curve of probability measures such that $m(0)=\mu_0$, $f,g,\psi$ be the corresponding functions in~\eqref{eq:fixed_funct}, and $(t,x)\mapsto u(t,x)$ be the value function~\eqref{eq:value_funct} for the optimization problem~\eqref{eq:control}--\eqref{eq:cost} or, equivalently, the unique bounded uniformly continuous viscosity solution of~\eqref{eq:hjb}. Recall that such $f$ satisfies~\eqref{eq:H'}. 

\n We consider the (local) measure--valued continuity equation
\begin{equation}\label{eq:transp}
\left\{
\begin{array}{l}
\partial_t \mu+\Div\big(\mu\,(-D_xu+f)\big)=0\,,\\
\mu(0)=\mu_0\,,
\end{array}
\right.
\end{equation}
where $\mu_0$ is any measure in $\prob_1(\R^d)$. Our aim is to construct a solution to~\eqref{eq:transp} via the method of characteristics, i.e. transporting the initial measure along integral curves of the vector field $v=-D_xu+f$ in~\eqref{eq:transp}, or equivalently along trajectories of the ODE~\eqref{eq:opt_control}.

As a preliminary, let us recall that, owing to Lemma~\eqref{lem:opt_set}--{\it (i)}, the multivalued map $(t,x)\mapsto\mathscr{A}(t,x)$ defined by~\eqref{eq:opt_set} has closed graph when $\Lsp^2([0,T]\,;\R^d)$ is endowed with the weak topology. Hence, we can find a Borel measurable selection $\frak{a}\colon[0,T]\times\R^d\to\Lsp^2([0,T]\,;\R^d)$ such that $\frak{a}(t,x)\in\mathscr{A}(t,x)$ for all $(t,x)$. 
We now define an absolutely continuous flow $\Phi(\cdot\,;t,x)$ by setting for all $s\in[t,T]$
\begin{equation}\label{eq:flow}
\Phi(s\,;t,x):= x+\int_t^s\left(\frak{a}(t,x)(\sigma)+f(\sigma,\Phi(\sigma\,;t,x))\right)\,d\sigma\,,
\end{equation}
and collect in next lemma the properties of $\Phi$ which we need to study~\eqref{eq:transp}.

\begin{lem}\label{lem:flow} Let $(t,x)\in[0,T[\,\times\R^d$, $u(t,x)$ be the semiconcave function given by~\eqref{eq:value_funct}, $f$ be of class $\con^2$ w.r.t. $x$ and satisfying~\eqref{eq:H'}, and $\Phi(\cdot\,;t,x)$ be the flow defined at~\eqref{eq:flow}. Then, the following facts hold.
\begin{description}
\item{\it (i)} For all $x\in\R^d$ and $t\leq s\leq s'\leq T$ there holds $\Phi(s'\,;t,x)=\Phi(s'\,;s,\Phi(s\,;t,x))$.
\item{\it (ii)} For all $x\in\R^d$ and $s\in(t,T)$ there holds $\partial_s \Phi(s\,;t,x) = - D_xu(s,\Phi(s\,;t,x))+ f(s,\Phi(s\,;t,x))$.
\item{\it (iii)} For all $x\in\R^d$ and $t\leq s\leq s'\leq T$ there holds
$$
\big|\Phi(s'\,;t,x)-\Phi(s\,;t,x)\big|\leq \big(\|D_xu\|_\infty+\|f\|_\infty\big) \,|s'-s|\,.
$$
\end{description}
Moreover, there exists a constant $K>0$, depending only on the constant $C$ in~\eqref{eq:H'}, such that
$$
|x-y|\leq K\,\big|\Phi(s\,;t,x)-\Phi(s\,;t,y)\big|\,,
\qquad\qquad
\mbox{for all }~x,y\in\R^d\,,~0\leq t<s\leq T\,,
$$
which implies, in particular, that $x\mapsto\Phi(s\,;t,x)$ has a Lipschitz continuous inverse on $\Phi(s\,;t,\R^d)$.
\end{lem}

\smallskip

\n{\it Proof.} Lemma~\ref{lem:opt_set}--{\it (iii)} implies that $\mathscr{A}(s,\Phi(s\,;t,x))=\{\frak{a}(t,x)_{|[s,T]}\}$ for all $t<s<T$. Thus, {\it (i)} follows.
In turn, Lemma~\ref{lem:opt_set}--{\it (ii)} implies that $u(s,\cdot)$ is differentiable at $\Phi(s\,;t,x)$ for $t<s<T$, and that $D_xu(s,\Phi(s\,;t,x))=-\frak{a}(t,x)(s)$, so that we conclude
$$
\partial_s \Phi(s\,;t,x) = \frak{a}(t,x)(s)+ f(s,\Phi(s\,;t,x))= - D_xu(s,\Phi(s\,;t,x))+ f(s,\Phi(s\,;t,x))\,,
$$
and {\it (ii)} is proved. Now, from {\it (ii)} we obtain
$$
\big|\Phi(s\,;t,x)-\Phi(s'\,;t,x)\big|\leq \int_s^{s'}\left|\partial_s \Phi(\sigma\,;t,x)\right|\,d\sigma
\leq\big(\|D_xu\|_\infty+\|f\|_\infty\big) \,|s'-s|\,
$$
which yields {\it (iii)}. Finally, $K$--Lipschitz continuity of the inverse can be proved as in Lemma~4.13 of~\cite{CardaNotes}, thanks to~\eqref{eq:H'} and to the semiconcavity of $u$ w.r.t. $x$.~~$\diamond$

\smallskip

Owing to Lemma~\ref{lem:flow} and to the fact that $\mu_0$ satisfies {\bf (H)}, we now claim that the curve of push--forward measures $s\mapsto\Phi(s\,;0,\cdot)\sharp \mu_0$ is the unique weak solution to~\eqref{eq:transp}. Recall that, given a Borel measurable function $h\colon\R^d\to\R^d$ and a measure $\mu\in\prob_1(\R^d)$, $h$ induces a probability measure $h\sharp\mu$ on $\R^d$, called \emph{push--forward measure} and defined by $h\sharp\mu(B):= \mu(h^{-1}(B))$, for all Borel subsets $B\subseteq \R^d$. The proof of the claim is the combination of the next two lemmas. 

\begin{lem}\label{lem:transp1} Let $\mu_0\in\prob_1(\R^d)$, $u$ be the semiconcave function given by~\eqref{eq:value_funct}, $f$ be of class $\con^2$ w.r.t. $x$ and satisfying~\eqref{eq:H'}, and $\Phi$ be the flow defined at~\eqref{eq:flow}. Then, $\mu(s):=\Phi(s\,;0,\cdot)\sharp \mu_0$ is a weak solution to~\eqref{eq:transp} and there holds 
\begin{equation}\label{eq:push_lip}
W_1(\mu(s),\mu(s'))\leq \big(\|D_xu\|_\infty+\|f\|_\infty\big) |s-s'|\,,\qquad\qquad\forall~t\leq s\leq s'\leq T\,.
\end{equation}
Moreover, if $\mu_0$ satisfies {\bf (H)}, i.e. $\mu_0\ll\Leb^d$ with density still denoted by $\mu_0$ which satisfies
$\mathrm{supp}~\mu_0\subseteq B(0,C)$ and
$\| \mu_0\|_\infty\leq C$, being $C$ the constant in {\bf (H)}, then $\mu(s)\ll\Leb^d$ for all $s\in[0,T]$, and there exists a constant $C_3>0$, depending on $C$ but not on $s$, such that the density functions, still denoted by $\mu(s)$, satisfy $\mathrm{supp}~\mu(s)\subseteq B(0,C_3)$ and $\| \mu(s)\|_\infty\leq C_3$ for all $s\in[0,T]$.
\end{lem}

\smallskip

\n{\it Proof.} The fact that $\Phi(s\,;0,\cdot)\sharp \mu_0$ is a weak solution to~\eqref{eq:transp} is standard (see e.g.~\cite{AmbrosioNotes}, or~\cite{CardaNotes}), so we focus on the other properties of such push--forward measure. By Lemma~\ref{lem:flow}--{\it (iii)}
$$
W_1(\mu(s),\mu(s'))\leq\int_{\R^d}|\Phi(s\,;0,x)-\Phi(s'\,;0,x)|\,d\mu_0(x)\leq\big(\|D_xu\|_\infty+\|f\|_\infty\big) |s-s'|\,.
$$
Recalling that $u$ is semiconcave w.r.t. $x$ with a constant $C_2$ only depending on $C$, we can exploit the bound $\|D_xu\|_\infty+\|f\|_\infty\leq C+ C_2$, together with the assumption $\mathrm{supp}~\mu_0\subseteq B(0,C)$, to conclude that all measures $\mu(s)$ with $s\in[0,T]$ satisfy
$\mathrm{supp}~\mu(s)\subseteq B(0,C_3)$,
for a suitable constant $C_3$ only depending on $C$ from {\bf (H)}.
Next, fixed $s\in[0,T]$, Lemma~\ref{lem:flow} implies that the map $x\mapsto\Phi(s\,;0,x)$ has a $K$--Lipschitz continuous inverse we denote with $\Xi$. Then, for every Borel subset $E\subseteq\R^d$ there holds
$$
\mu(s)(E)=\mu_0\big(\Phi(s\,;0,\cdot)^{-1}(E)\big)=\mu_0(\Xi(E))\leq\|\mu_0\|_\infty\Leb^d(\Xi(E))\leq\|\mu_0\|_\infty\,K^d\Leb^d(E)\,,
$$
whence both $\mu(s)\ll\Leb^d$ and $\|\mu(s)\|_\infty\leq K^d\,\|\mu_0\|_\infty$ for all $s\in[0,T]$ follow.~~$\diamond$

\smallskip

\begin{lem}\label{lem:transp2} Assume that $\mu_0$ satisfies {\bf (H)} and let $u$ be the semiconcave function given by~\eqref{eq:value_funct}, $f$ be of class $\con^2$ w.r.t. $x$ and satisfying~\eqref{eq:H'}, and $\Phi(\cdot\,;t,x)$ be the flow defined at~\eqref{eq:flow}.
Then, $\mu(s):=\Phi(s\,;0,\cdot)\sharp \mu_0$ is the unique weak solution to~\eqref{eq:transp}.
\end{lem}

\smallskip

\n{\it Proof.} We claim that every solution ${\frak m}(\cdot)\in\con([0,T],\prob_1(\R^d))$ of the continuity equation~\eqref{eq:transp} is a ``superposition solution'', i.e., it admits the following probabilistic representation
\begin{equation}\label{eq:eta_probab}
\int_{\R^d}\vp(x)\,d{\frak m}(t)(x)=\int_{\R^d\times \Gamma_T}\vp(\gamma(t))\,d\eta(x,\gamma)\,,
\qquad\qquad
\forall~\vp\in\con_b(\R^d)\,,
\end{equation}
where $\Gamma_T:= \con([0,T],\R^d)$, $\con_b(\R^d)$ is the set of bounded continuous functions from $\R^d$ to $\R$, and $\eta\in\prob(\R^d\times \Gamma_T)$ is a probability measure concentrated on Carath\'eodory solutions $\gamma(\cdot)$ of
\begin{equation}\label{eq:Cauchy_help}
 \dot \gamma(s) = -D_xu(s,\gamma(s))+f(s,\gamma(s))\,,
 \qquad\qquad
 \gamma(0)=x\,.
\end{equation}
Since the vector field in~\eqref{eq:Cauchy_help} is bounded, the proof of such claim is a direct application of the so called \emph{superposition principle} for measure--valued solutions to~\eqref{eq:cont_eq} (cf Theorem~12 of~\cite{AmbrosioNotes}).
Then, observe that the first marginal of such a measure $\eta$ is $\mu_0$, because for all $\vp\in\con_b(\R^d)$
$$
\int_{\R^d\times \Gamma_T}\vp(x)\,d\eta(x,\gamma)=\int_{\R^d\times \Gamma_T}\vp(\gamma(0))\,d\eta(x,\gamma)=\int_{\R^d}\vp(x)\,d{\frak m}(0)(x)=\int_{\R^d}\vp\,d\mu_0\,.
$$
Hence, we can disintegrate the measure $\eta$ w.r.t. $\mu_0$, and find a $\mu_0$--a.e. uniquely determined family of measures $(\lambda_x)$ on $\Gamma_T$ such that for all Borel measurable maps $\zeta\colon\R^d\times\Gamma_T\to[0,+\infty]$ there holds
\begin{equation}\label{eq:eta_disint}
\int_{\R^d\times\Gamma_T} \zeta(x,\gamma)\,d\eta(x,\gamma)=\int_{\R^d}\!\int_{\Gamma_T}\zeta(x,\gamma)\,d\lambda_x(\gamma)\,d\mu_0(x)\,.
\end{equation}
In particular, $\lambda_x$ is concentrated on trajectories of~\eqref{eq:Cauchy_help}. Since Proposition~\ref{prop:iff_optimal_control} implies that such trajectory is unique for $\mu_0$--a.e. $x\in\R^d$, 
because ${\cal A}(0,x)$ is a singleton ($\Leb^d$-- and thus) $\mu_0$--almost everywhere, we can conclude that the measure $\lambda_x$ is a Dirac delta concentrated on the unique solution $\Phi(\cdot\,;0,x)$ for $\mu_0$--a.e. initial data $x\in\R^d$.
This implies that $\int_{\R^d}\vp(x)\,d{\frak m}(t)(x)=\int_{\R^d}\vp(\Phi(t\,;0,x))\,d\mu_0(x)$, whence the uniqueness of the solution to~\eqref{eq:transp} follows.~~$\diamond$

\subsection{The fixed point argument}\label{sec:PROOF}

We finally define the function ${\cal T}$ whose fixed points are solutions to~\eqref{eq:mfg}. Let ${\cal K}$ be the convex subset of $\con([0,T],\prob_1(\R^d))$ defined by
\begin{align*}
{\cal K}:=\big\{m(\cdot)&\in\con([0,T],\prob_1(\R^d))~;~m(0)=\mu_0 \mbox{ and $m(\cdot)$ is Lipschitz with }Lip(m)\leq C+C_2\big\}\,,
\end{align*}
where $C$ is the constant in {\bf (H)} and $C_2$ is the constant in~\eqref{eq:semiconc}, and let us fix $m(\cdot)\in{\cal K}$. 
By applying the results from the previous sections, the HJB equation
$$
-\partial_t u(t,x)+\,{1\over 2}\,|D_xu(t,x)|^2-D_xu(t,x)\cdot F(x,m(t))-G(x,m(t))=0\,,
$$
with final data $u(T,x)=\Psi(x,m(T))$, admits a unique Lipschitz continuous viscosity solution $\bar u$ and such a solution is also semiconcave in the variable $x$, with a semiconcavity constant $C_2$ independent from the curve $m(\cdot)$. By applying Lemma~\ref{lem:transp2}, we also deduce that the continuity equation
$$
\partial_t \mu(t)+\Div\big(\mu(t)\,(-D_x\bar u(t,x)+F(x,m(t))\big)=0\,,\qquad\qquad \mu(0)=\mu_0\,,
$$
has a unique solution $\bar\mu$ and that such a solution belongs to ${\cal K}$ thanks to Lemma~\ref{lem:transp1}.
Hence, it is well--defined a map ${\cal T}\colon{\cal K}\to{\cal K}$ by setting $m(\cdot)\mapsto {\cal T}(m(\cdot)):= \bar\mu(\cdot)$. Moreover, ${\cal T}$ is compact, because the curve $s\mapsto\bar\mu(s)$ is uniformly Lipschitz continuous with values in the compact set $\prob_1(B(0,C_3))$, where $C_3>0$ is the constant introduced in Lemma~\ref{lem:transp1}.

\n If we prove that ${\cal T}$ is continuous, w.r.t. the uniform topology on $\con([0,T],\prob_1(\R^d))$, then we can conclude by applying Schauder's fixed point theorem. Thus, we consider a sequence $m_n(\cdot)\in{\cal K}$ which converges uniformly to $m(\cdot)\in\con([0,T],\prob_1(\R^d))$. Since the curves $m_n(\cdot)$ have Lipschitz constant uniformly bounded, the limit $m(\cdot)$ is still an element of ${\cal K}$. Moreover, for each $n\in\N$, let us denote with $\bar u_n$ the sequence of uniformly semiconcave solutions to the HJB equation
$$
-\partial_t u(t,x)+\,{1\over 2}\,|D_xu(t,x)|^2-D_xu(t,x)\cdot F(x,m_n(t))-G(x,m_n(t))=0\,,
$$
with final data $u(T,x)=\Psi(x,m_n(T))$, 
and with $\bar \mu_n$ the sequence of solutions to
$$
\partial_t \mu(t)+\Div\big(\mu(t)\,(-D_x \bar u_n(t,x)+F(x,m_n(t))\big)=0\,,\qquad\qquad \mu(0)=\mu_0\,.
$$
We also denote with $\bar u$ the solution to HJB corresponding to the limit curve $m(\cdot)$ and with $\bar \mu$ the solution to the continuity equation corresponding to $\bar u$ and $m(\cdot)$. 

We claim that $\bar u_n\to \bar u$ locally uniformly on $[0,T]\times\R^d$ and that $\bar \mu_n\to\bar \mu$ in $\con([0,T],\prob_1(\R^d))$, whence continuity of ${\cal T}$ follows.
Indeed, assumptions {\bf (H)} on $F,G,\Psi$ ensure that the sequences of maps $F(x, m_n(t))$, $G(x, m_n(t))$ and $\Psi(x, m_n(T))$ locally uniformly converge to the maps $F(x, m(t))$, $G(x, m(t))$ and $\Psi(x, m(T))$. Hence, the convergence of the solutions $\bar u_n(t,x)$ to $\bar u(t,x)$ locally uniformly in $[0,T]\times\R^d$ follows by standard results on viscosity solutions (cf~\cite{CS}). Moreover, each $\bar u_n$ is semiconcave wr.t. $(t,x)$ with a common constant dependent only on $C$ and $C+C_2$. Therefore, we deduce by Proposition~\ref{prop:semiconc_lip} that $D_x\bar u_n\to D_x\bar u$ $\Leb^{d+1}$--a.e. in $[0,T]\times\R^d$.

\n Since $\mu_0$ is assumed to be absolutely continuous w.r.t. $\Leb^d$, all curves $\bar \mu_n$ consist of absolutely continuous measures and they satisfy the properties in Lemma~\ref{lem:transp1}. By applying Ascoli--Arzel\`a theorem and Banach--Alaoglu theorem, we can extract a subsequence, still denoted with $\bar \mu_n$, which converges in $\con([0,T],\prob_1(B(0,C_3)))$ to a certain curve of measures $\hat\mu(\cdot)$, and whose densities converge in $w^*-\Lsp^\infty([0,T]\times\R^d)$ to the densities of the measures $\hat\mu(\cdot)$.

\n Since the curves $\bar \mu_n$ solve the continuity equation for $\bar u_n$ and $m_n$, and for any test function $\vp$
$$
\partial_s \vp+\langle D_x \vp,-D_x\bar u_n+F(\cdot,m_n)\rangle~\overset{\Lsp^1((0,T)\times\R^d)}{\longrightarrow}~\partial_s \vp+\langle D_x \vp,-D_x\bar u+F(\cdot,m)\rangle
$$
as $n\to+\infty$, one concludes that $\hat\mu$ is a weak solution to the continuity equation for $\bar u$ and $m$. However, the solution to this limit equation is unique by Lemma~\ref{lem:transp2}, and such uniqueness implies that ${\cal T}(m_n(\cdot))\to\hat\mu(\cdot)=\bar \mu(\cdot)={\cal T}(m(\cdot))$ as $n\to\infty$.~~$\diamond$

\section{Extensions and open problems}\label{sec:fine}

\subsection{Multipopulation models}

Our results extend easily to the case of several different populations of pedestrians. Given its interest in applications, we give the basic details here. Assume that $M$ populations are present in the walking area, with continuous distributions $\varrho^{(1)},\ldots,\varrho^{(M)}$, and that they obey continuity equations 
$$
\partial_t\varrho^{(j)}+\mathrm{div}(\varrho^{(j)} v^{(j)}) = 0\,,
\qquad
\varrho^{(j)}(0,\cdot)=\varrho^{(j)}_0\,,
\qquad\qquad j=1,\ldots,M\,,
$$
with suitable initial distributions $\varrho^{(1)}_0,\ldots,\varrho^{(M)}_0$ and with velocities $v^{(j)}$ of the same form as~\eqref{eq:v}. To fix the ideas, let us assume that such velocities are explicitly given as superposition of a behavioral velocity and several interaction velocities accounting for the repulsion/attraction between individuals of the $j$th population with individuals of the others
\begin{align*}
v^{(j)}(t,x)~&:=~v_b^{(j)}(t,x) ~+~\sum_{k=1}^Mv^{(j,k)}_i[\rho^{(k)}(t)](x)\,,
\end{align*}
for instance with integral interaction velocities
\begin{align*}
v^{(j,k)}_i(t,x)~&= ~ \int_{B(x,R)}\F_{jk}(y-x)\,\varrho^{(k)}(t,y)\,dy\,,
\end{align*}
with functions $\F_{jk}$ as in~\eqref{eq:vi}. Here, as in Section~\ref{sec:theory_n_pedestrians}, we used the notations $\rho^{(j)}(t):=\,{\varrho^{(j)}(t,\cdot)\over\|\varrho^{(j)}_0\|_{\Lsp^1(\R^2)}}\,\Leb^2$. 
We can still assume that pedestrians are \emph{highly rational} when they have to choose their behavioral velocities $v_b^{(j)}(t,x)$, i.e. that they are able to exploit the full knowledge of all distribution measures $\rho^{(1)},\ldots,\rho^{(M)}$ in $[0,T]$. Thus, if we fix some performance criteria of the form
$$
{\cal C}^{(j)}(t,x, v_b(\cdot))~:=~\int_t^T \left({1\over 2} |v_b(s)|^2+\ell^{(j)}[\rho^{(j)}(s)](y^{(j)}(s))\right)\,ds+h^{(j)}[\rho^{(j)}(T)](y^{(j)}(T))\,,
$$
with cost functions $\ell^{(j)}$ and $h^{(j)}$ possibly different among the populations but all satisfying {\bf (H2)}, and with $y^{(j)}(\cdot)$ being the Carath\'eodory solution to
$$
\dot y^{(j)}(t)=v^{(j)}(t,y^{(j)}(t))\,,
$$
we can repeat the construction outlined in the previous sections. Namely, pedestrians of each population define their own value function as $\phi^{(j)}(t,x)\!:=\inf_{v_b(\cdot)}{\cal C}^{(j)}(t,x, v_b(\cdot))$ and they can find it by solving an associated HJB equation with Hamiltonian
$$
H^{(j)}(t,x,D_x\phi^{(j)}):=\max_{v_b \in \R^2}\left\{ -(v_b+v_i^{(j)}(x)\big)\cdot D_x\phi^{(j)}(t,x)-{1\over 2} |v_b|^2-\ell^{(j)}(x)\right\}\,,
$$
with final condition $\phi^{(j)}(T,x)=h^{(j)}(x)$, where we have skipped the non--local dependences in $v_i^{(j)}$, $\ell^{(j)}$ and $h^{(j)}$ for sake of brevity. Once $\phi^{(j)}$ is known, the behavioral field $v_b^{(j)}$ takes again the form $v_b^{(j)}=-D_x\phi^{(j)}$, so that the new mean field system for the multipopulation problem in $(0,T)\times\R^2$ becomes the following system of $2M$ forward--backward equations, for $j=1,\ldots,M$, 
\begin{equation}\label{eq:mfg3}
\left\{
\begin{array}{l}
-\partial_t\phi^{(j)}(t,x)+\displaystyle\,{1\over 2}\,|D_x\phi^{(j)}(t,x)|^2-\sum_{k=1}^Mv^{(j,k)}_i[\rho^{(k)}(t)](x)\cdot D_x\phi^{(j)}(t,x)=\ell[\rho^{(j)}(t)](x)\\
\partial_t \rho^{(j)}(t)+\Div\Bigg(\rho^{(j)}(t)\bigg(-D_x\phi^{(j)}(t,x)+ \displaystyle\sum_{k=1}^Mv^{(j,k)}_i[\rho^{(k)}(t)](x)\bigg)\Bigg)=0\phantom{\displaystyle{1\over 2}}\\
\rho^{(j)}(0)=\varrho^{(j)}_0/\|\varrho^{(j)}_0\|_{\Lsp^1}\,\Leb^2\,,~~~\phi^{(j)}(T,\cdot)=h^{(j)}[\rho^{(j)}(T)]\phantom{\displaystyle{1\over 2}}
\end{array}
\right.
\end{equation}
In this case, our main theorem rewrites as follows.

\begin{thm} Assume that for, all $j,k=1,\ldots,M$, the interaction velocity fields $v^{(j,k)}_i$ satisfy {\bf (H1)}, the cost functions $\ell^{(j)}$ and $h^{(j)}$ satisfy {\bf (H2)}, and the initial distributions $\varrho^{(j)}_0$ belong to $\Lsp^1(\R^2)\cap\Lsp^\infty(\R^2)$ and have compact support. Then, system
\eqref{eq:mfg3} admits a solution $2M$--tuple $(\phi^{(j)},\varrho^{(j)})\in W^{1,\infty}((0,T)\times\R^2)\times\Lsp^1((0,T)\times\R^2)$, with the functions $\phi^{(j)}$ solving the HJB backward equations in the viscosity sense and with the functions $\varrho^{(j)}$ solving the forward continuity equations in the sense of distributions.
\end{thm}

\n The proof is based on a simple extension of Theorem~\ref{thm:mfg} to the case of curves in $\con\big([0,T]\,;(\prob_1(\R^2))^M\big)$ and therefore is omitted. Further extensions to more general couplings between the populations
$$
v_i^{(j)}=v_i^{(j)}[\rho^{(1)},\ldots,\rho^{(M)}]\,,
\quad 
\ell^{(j)}=\ell^{(j)}[\rho^{(1)},\ldots,\rho^{(M)}]\,,
\quad 
h^{(j)}=h^{(j)}[\rho^{(1)},\ldots,\rho^{(M)}]\,,
$$ 
can be proved analogously, as long as assumptions {\bf (H1)} and {\bf (H2)} are replaced with continuity of the applications as maps defined on $[\prob_1(\R^2)]^M$ and the bounds on their $\con^2$ norm are uniform w.r.t. any $M$--tuple of probability measures.

\subsection{Open problems} 

In this paper we have showed that the general mean field system~\eqref{eq:mfg} admits a solution pair $(u,\mu)$ defined in $(0,T)\times\R^d$, and we have seen how this applies to a model for $2$d crowd dynamics with macroscopic velocity~\eqref{eq:v}, for a suitable choice of the non--local interaction field $v_i$ in~\eqref{eq:vi}. However, this is only a first step in the analysis of what was presented as a model for highly rational pedestrians in~\cite{CrPrTo}, and many questions remain unanswered. 

First of all, we do not know yet whether the solution to~\eqref{eq:mfg2} is unique or if the possibility to take very large desired velocities $v_b^*$ allows for the presence of multiple Nash equilibria strategies, and thus of multiple solution pairs $(\rho,\phi)$. The coupling between the equations in~\eqref{eq:mfg} makes impossible to directly apply the monotonicity arguments by~\cite{LL1}, and at present it is not clear which stronger assumptions on $v_i,\ell,g$ might lead to uniqueness without losing physical relevance. 

A second challenge is the handling of bounded domains. In Theorem~\ref{thm:crowd} we were only capable to cover the case of dynamics taking place in the whole $\R^2$, even if a posteriori the support of the solution $\rho$ remains in a suitably large compact set, whenever the initial distribution is compactly supported. When trying to adapt the results to the case of a bounded domain $\Omega$, as a room or a railway station, we have to face several difficulties. First of all, the handling of the boundary conditions: even if some results are available for both HJB equations~\cite{Soner1} and continuity equations~\cite{CrippaSpinolo}, their actual combination for system~\eqref{eq:mfg2} is to our knowledge an open problem, outside the case of periodic data on a torus, which is not very meaningful for applications to crowd dynamics. Also, when one wants to describe the case of a bounded environment $\Omega$ with a target $\Sigma\subset\partial\Omega$ (as is typically the case when modeling evacuation of people from a building or from a room), one has to combine some viability condition on $\partial\Omega\setminus \Sigma$, so to model ``impermeability'' of the walls of the environment, with suitable outflow conditions on $\Sigma$, so to model the actual evacuation, and this typically reduces the regularity of the optimal vector field $v_b^*$, making harder to solve the continuity equation.

Finally, in view of applications to real crowd dynamics, it would be useful to extend the theory to problems where the performance criterion is the minimal evacuation time, and problems where the behavioral velocity is subject to additional restrictions (e.g. it can only attains values in a fixed bounded subset ${\bf U}$ of $\R^2$).


\begin{thebibliography}{99}


\bibitem{AmbrosioNotes} L. Ambrosio, \emph{Transport equation and Cauchy problem for non-smooth vector fields}. in B. Dacorogna, P. Marcellini (Eds.), Calculus of variations and nonlinear partial differential equations, 1--41. Springer, Berlin, 2008.

\bibitem{bellomo2011} N. Bellomo \& C. Dogb\'{e}, \emph{On the modeling of traffic and crowds: A survey of models, speculations, and perspectives}, SIAM Rev. {\bf 53}, n. 3, 409--463, 2011.

\bibitem{BuDFMaWo} M. Burger, M. Di Francesco, P. A. Markowich \& M.-T. Wolfram, \emph{Mean field games with nonlinear mobilities in pedestrian dynamics}, submitted, 2013.

\bibitem{BuFr} G. Buttazzo \& L. Freddi, \emph{Optimal control problems with weakly converging input operators}, Discr. Cont. Dyn. Syst. {\bf 1}, n. 3, 401--420, 1995.

\bibitem{CS} P. Cannarsa \& C. Sinestrari, Semiconcave Functions, Hamilton--Jacobi Equations, and Optimal Control. 
Birkh\"auser Boston, Inc., Boston, MA, 2004.

\bibitem{CardaNotes} P.~Cardaliaguet, Notes on Mean Field Games (from P.-L. Lions' lectures at Coll\`ege de France), 2010.

\bibitem{CHM} R.~M.~Colombo, M.~Herty \& M.~Mercier, \emph{Control of the continuity equation with a non local flow},  ESAIM Control Optim. Calc. Var. {\bf 17}, 353--379, 2011

\bibitem{CrippaSpinolo} G.~Crippa, C.~Donadello \& L.~V.~Spinolo, \emph{Initial--boundary value problems for continuity equations with $BV$ coefficients},  J. Math. Pures Appl. To appear.

\bibitem{CPT10}E. Cristiani, B. Piccoli \& A. Tosin, {\it Modeling self-organization in pedestrians and animal groups from macroscopic and microscopic viewpoints}, in G. Naldi, L. Pareschi, G. Toscani (Eds.), ÒMathematical modeling of collective behavior in socio-economic and life sciencesÓ, 337-364, BirkhŠuser Boston, Cambridge, MA, 2010.

\bibitem{CPT11} E. Cristiani, B. Piccoli \& A. Tosin,
{\it Multiscale modeling of granular flows with application to crowd dynamics}, Multiscale Model. Simul., {\bf 9}, n. 1, 155--182, 2011.

\bibitem{cristiani2012CDC} E. Cristiani, B. Piccoli \& A. Tosin,
{\it How can macroscopic models reveal self-organization in traffic flow?}, in Proceedings of the $51$st IEEE Conference on Decision and Control, Maui, HI, USA, 6989--6994, 2012.

\bibitem{CPTbook} E. Cristiani, B. Piccoli \& A. Tosin, Multiscale Modeling of Pedestrian Dynamics. Springer, Berlin Heidelberg, 2014.

\bibitem{CrPrTo} E. Cristiani, F. S. Priuli \& A. Tosin, {\it Modeling rationality to control self-organization of crowds: an environmental approach}, preprint, 2014.

\bibitem{duives2013} D. C. Duives, W. Daamen \& S. P. Hoogendoorn, {\it State--of--the--art crowd motion simulation models}, Transportation Res. C, in press, 2013.

\bibitem{Dogbe} C. Dogb\'{e}, \emph{Modeling crowd dynamics by the mean-field limit approach}, Math. Comput. Modelling {\bf 52}, n. 9--10, 1506--1520, 2010.

\bibitem{GLL} O.~Gu\'eant, J.-M.~Lasry, P.-L.~Lions,  \emph{Mean field games and applications}, in ``Paris-Princeton Lectures on Mathematical Finance 2010'', R.A. Carmona et al. eds.,
205--266, Springer, Berlin Heidelberg, 2011.

\bibitem{helbing2001} D. Helbing, \emph{Traffic and related self-driven many-particle systems}, Rev. Mod. Phys. {bf 73}, n. 4, 1067--1141, 2001.

\bibitem{HB} S. P. Hoogendoorn \& P. H. L. Bovy, {\it Simulation of pedestrian flows by optimal control and differential games}, Optimal Control Appl. Methods, {\bf 24}, 153--172, 2003.

\bibitem{HB2} S. P. Hoogendoorn \& P. H. L. Bovy, {\it Dynamic user--optimal assignment in continuous time and space}, Transportation Res. B, {\bf 38}, 571--592, 2004.

\bibitem{Hughes}  R. L. Hughes, {\it A continuum theory for the flow of pedestrians}, Transportation Res. B, {\bf 36}, 507--535, 2002.

\bibitem{kachroo2008} P. Kachroo,  S. J. Al-nasur, S. A. Wadoo \& A. Shende, Pedestrian dynamics. Feedback control of crowd evacuation. 
Springer, Berlin Heidelberg, 2011

\bibitem{LW} A. Lachapelle \& M.-T. Wolfram, \emph{On a mean field game approach modeling congestion and aversion in pedestrian crowds},
Transportation Res. B, {\bf 45}, n. 10, 1572--1589, 2011.

\bibitem{LL1} J.-M.~Lasry \& P.-L.~Lions, \emph{Jeux \`a champ moyen. I. Le cas stationnaire.} C. R. Acad. Sci. Paris {\bf 343}, 619--625, 2006.

\bibitem{LL2} J.-M.~Lasry \& P.-L.~Lions, \emph{Jeux \`a champ moyen. II. Horizon fini et contr\^ole optimal.} C. R. Acad. Sci. Paris {\bf 343}, 679--684, 2006.

\bibitem{PiccoliRossi} B.~Piccoli \& F.~Rossi, \emph{Transport equation with nonlocal velocity in Wasserstein spaces:
convergence of numerical schemes}, Acta Appl. Math. {\bf 124}, 73--105, 2013.

\bibitem{Soner1} H.~M.~Soner, \emph{Optimal Control with state--space constraint I}, SIAM J. Cont. Opt. {\bf 24}, n. 3, 552--562, 1986.

\bibitem{TosinFrasca} A.~Tosin \& P.~Frasca, \emph{Existence and approximation of probability measure solutions to models of collective behaviors},
Netw. Heterog. Media {\bf 6}, n. 3, 561--596, 2011.


\end{thebibliography}
\end{document}